\documentclass[11pt,reqno]{amsart}

\usepackage{amsmath,amsthm}
\usepackage{amssymb}
\usepackage{euscript}
\usepackage[frame,ps,matrix,arrow,curve,rotate]{xy}

\numberwithin{equation}{subsection}

\newcommand{\sqsp}{\renewcommand{\baselinestretch}{1.15}\tiny\normalsize}

\newcommand{\relphantom[1]}{\mathrel{\phantom{#1}}}
\newtheorem{thm}[subsection]{Theorem}
\newtheorem{lemma}[subsection]{Lemma}

\newtheorem{cor}[subsection]{Corollary}

\theoremstyle{definition}

\newcommand{\biglbrack}{\Biggl \lbrack}  % big [
\newcommand{\bigrbrack}{\Biggr \rbrack}  % big ]
\newcommand{\biglpren}{\bigl (}  % medium big (
\newcommand{\bigrpren}{\bigr )}  % medium big )
\newcommand{\Mbar}{\overline{M}}
\newcommand{\Mtilde}{\widetilde{M}}
\newcommand{\mbar}{\overline{m}}

\newcommand{\fbar}{\overline{f}}

\newcommand{\Fbar}{\overline{F}}
\newcommand{\Ftilde}{\widetilde{F}}
\newcommand{\Thetabar}{\overline{\Theta}}
\newcommand{\Thetatilde}{\widetilde{\Theta}}
\newcommand{\thetabar}{\overline{\theta}}

\newcommand{\sumprime}{\sideset{}{'}\sum}  % \sum^\prime
\newcommand{\Hzinb}{H_{\text{Zinb}}}
\newcommand{\Czinb}{C_{\text{Zinb}}}

\DeclareMathOperator{\Id}{Id}
\DeclareMathOperator{\Hom}{Hom}

\DeclareMathOperator{\Ob}{Ob}

\begin{document}
\title{Deformation of dual Leibniz algebra morphisms}
\author{Donald Yau}

\begin{abstract}
An algebraic deformation theory of morphisms of dual Leibniz algebras is obtained.
\end{abstract}

%\subjclass[2000]{}
\email{dyau@math.ohio-state.edu}
\address{Department of Mathematics, The Ohio State University at Newark, 1179 University Drive, Newark, OH 43055, USA}

%\date{\today}
\maketitle
\sqsp

%%================%%
%%                %%
%%  Introduction  %%
%%                %%
%%================%%

\section{Introduction}
\label{sec:intro}

A Leibniz algebra, as introduced by Loday \cite{loday1,loday}, is a vector space $R$ together with a bilinear binary operation $\lbrack -, - \rbrack$ for which $\lbrack -, z \rbrack$ is a derivation for every $z \in R$.  Thus, Leibniz algebras are non-commutative versions of Lie algebras, which are Leibniz algebras whose brackets are skew-symmetric.  Algebraic deformations, in the sense of Gerstenhaber \cite{ger1}, of Leibniz algebras have been studied by Balavoine \cite{bal1,bal2,bal3}.  In fact, the results in \cite{bal3} applies to algebras over any binary quadratic operad.  Dual Leibniz algebras \cite{loday0,loday}, also known as \emph{Zinbiel algebras}, are the operadic dual \cite{gk} of Leibniz algebras, and the operad that defines Zinbiel algebras is binary quadratic.  In particular, Balavoine's operadic approach to algebraic deformations applies to Zinbiel algebras.

The purpose of this paper is to study algebraic deformations of \emph{morphisms} of Zinbiel algebras.  In doing so, we use the results of Balavoine \cite{bal3} freely.  Our results are analogous to those of Gerstenhaber and Schack \cite{gs1,gs2,gs3}, who worked with morphisms and, more generally, diagrams of associative algebras.  We note that deformations of Lie algebra morphisms have been studied by Nijenhuis and Richardson \cite{nr} and, more recently, by Fr\'eiger \cite{freiger}.

The rest of this paper is organized as follows.  In Section \ref{sec:def complex}, the deformation complex $\Czinb^*(f,f)$ of a morphism $f$ of dual Leibniz algebras is introduced.  Deformations and their infinitesimals of a morphism $f$ are introduced in Section \ref{sec:def inf}.  Given a deformation, its infinitesimal is identified with a $2$-cocycle in the deformation complex, whose cohomology class is determined by the equivalence class of the deformation (Theorem \ref{thm:inf}).

A morphism $f$ is said to be \emph{rigid} if every deformation of it is equivalent to the trivial deformation.  In Section \ref{sec:rigidity}, a cohomological criterion for rigidity is given.  More precisely, it is shown that the vanishing of $\Hzinb^2(f,f) = H^2(\Czinb^*(f,f))$ implies that $f$ is rigid (Theorem \ref{thm:rigidity} and Corollary \ref{cor:rigidity}).  Since every infinitesimal is a $2$-cocycle, it is natural to ask, given a $2$-cocycle, under what condition is it the infinitesimal of a deformation.  This question is considered in Section \ref{sec:extending} and Section \ref{sec:Lemma}.  It is shown that there is an associated sequence of obstruction $3$-cocycles (Lemma \ref{lem:extending}).  The simultaneous vanishing of their cohomology classes is equivalent to the existence of a deformation with the given $2$-cocycle as its infinitesimal (Theorem \ref{thm:extending}).  In particular, the existence of such a deformation is guaranteed if $\Hzinb^3(f,f)$ is trivial (Corollary \ref{cor:extending}).

The results in \cite{bal3} will be used below, often implicitly, whenever we are dealing with a deformation of an individual Zinbiel algebra.

%%=======================%%
%%                       %%
%%  Deformation Complex  %%
%%                       %%
%%=======================%%

\section{Deformation complex}
\label{sec:def complex}

Before we introduce the deformation complex of a morphism, we briefly review the relevant definitions from \cite{loday0,loday} concerning dual Leibniz algebras and their cohomology.

%%%%%%%%%%%%%%%%%%%%%%%%%%%%%%%%%
\subsection{Dual Leibniz algebra}
\label{subsec:dual}

Throughout this paper we work with an arbitrary but fixed field $K$.  A \emph{dual Leibniz algebra} \cite{loday0,loday}, or a \emph{Zinbiel algebra}, is a $K$-vector space $R$ which comes equipped with a $K$-bilinear binary operation $R^{\otimes 2} \to R$, $(x,y) \mapsto x \cdot y$, satisfying the dual Leibniz identity
   \begin{equation}
   \label{eq:Leibniz}
   (x \cdot y) \cdot z ~=~ x \cdot (y \cdot z) \,+\, x \cdot (z \cdot y)
   \end{equation}
for $x, y, z \in R$.  Sometimes we write $m_R(x,y)$ for the product $x \cdot y$ and $(R,m_R)$ for the dual Leibniz algebra if there are at least two dual Leibniz algebras under consideration.

If $R$ and $S$ are Zinbiel algebras, then a morphism $f \colon R \to S$ is a $K$-linear map of the underlying vector spaces that respects the products, in the sense that $f(x \cdot y) = f(x) \cdot f(y)$ for $x, y \in R$.

Let $R$ be a Zinbiel algebra.  An \emph{$R$-bimodule} is a $K$-vector space $A$ together with $K$-bilinear operations $R \otimes A \to A$, $(r,a) \mapsto r \cdot a$, and $A \otimes R \to A$, $(a,r) \mapsto a \cdot r$, such that \eqref{eq:Leibniz} holds whenever exactly one of $x, y, z$ is from $A$ and the others from $R$.  For example, if $g \colon R \to S$ is a morphism of Zinbiel algebras, then $S$ can be given an $R$-bimodule structure, where $R \otimes S \to S$ is defined by $(r,s) \mapsto r \cdot s = g(r) \cdot s$ and similarly for $S \otimes R \to S$.  In particular, when $g$ is the identity map on $R$, this gives an $R$-bimodule structure on $R$.

%%%%%%%%%%%%%%%%%%%%%%%%%%%%%%%%%%%%%%%
\subsection{Low dimensional cohomology}
\label{subsec:coh}

Let $R$ be a Zinbiel algebra and let $A$ be an $R$-bimodule.  For the purpose of studying algebraic deformations, only low dimensional cohomology is needed, which we now recall.  Define a cochain complex $(\Czinb^*(R,A), d)$, $1 \leq * \leq 4$, as follows.  For $1 \leq n \leq 4$, set
   \[
   \Czinb^n(R,A) = \Hom_K(R^{\otimes n}, A).
   \]
For $1 \leq i \leq 3$, define the maps
   \[
   d^i \colon \Czinb^i(R,A) \to \Czinb^{i+1}(R,A)
   \]
by:
   \[
   \allowdisplaybreaks
   \begin{split}
   (d^1\varphi)(x, y) &~=~ x \cdot \varphi(y) \,-\, \varphi(x \cdot y) \,+\, \varphi(x)\cdot y,    \\
   (d^2\varphi)(x, y, z) &~=~ x \cdot (\varphi(y,z) \,+\, \varphi(z,y)) \,-\, \varphi(x\cdot y, z) \\
                   &\relphantom{} \relphantom{} ~+~ \varphi(x, y\cdot z + z \cdot y) \,-\, \varphi(x, y) \cdot z, \\
   (d^3\varphi)(x, y, z, w) &~=~ x \cdot \lbrace \varphi(y, z, w) \,-\, \varphi(z, w, y) \,+\, \varphi(z, y, w) \,-\, \varphi(w, z, y) \rbrace \\
                      &\relphantom{} \relphantom{} ~-~ \varphi(x\cdot y, z, w) \,+\, \varphi(x, y\cdot z + z \cdot y, w) \\
                      &\relphantom{} \relphantom{} ~-~ \varphi(x, y, z \cdot w + w \cdot z) \,+\, \varphi(x, y, z) \cdot w.
   \end{split}
   \]
It is easy to check directly that $d^{i+1} d^i = 0$ for $i = 1, 2$, and so $(\Czinb^*(R,A), d)$, $1 \leq * \leq 4$, is a cochain complex, which is usually abbreviated to $\Czinb^*(R,A)$.  For $n = 2, 3$, the cohomology group $\Hzinb^n(R,A)$ of $R$ with coefficients in the bimodule $A$ is defined to be $H^n(\Czinb^*(R,A))$.

%%%%%%%%%%%%%%%%%%%%%%%%%%%%%%%%%%%%%%%%%%%%%%
\subsection{Deformation complex of a morphism}
\label{subsec:def morphism}

Throughout the rest of this paper, $f \colon R \to S$ denotes an arbitrary but fixed morphism of Zinbiel algebras.  We consider $S$ also as an $R$-bimodule via $f$ wherever appropriate.

Setting $\Czinb^0(R,S) = 0$, we define the vector spaces $\Czinb^n(f,f)$ for $1 \leq n \leq 4$ as
   \begin{equation}
   \label{eq:def complex}
   \Czinb^n(f,f) ~\buildrel \text{def} \over =~ \Czinb^n(R,R) \times \Czinb^n(S,S) \times \Czinb^{n-1}(R,S).
   \end{equation}
For $1 \leq i \leq 3$, the differential $d^i_f \colon \Czinb^i(f,f) \to \Czinb^{i+1}(f,f)$ is defined by
   \begin{equation}
   \label{eq:d}
   d^i_f(\xi; \pi; \varphi) ~\buildrel \text{def} \over =~ (d^i \xi; d^i\pi; f \xi - \pi f - d^{i-1} \varphi).
   \end{equation}
Here $f\xi, \pi f \in \Czinb^i(f,f)$ are the push-forwards given by
   \[
   \begin{split}
   (f\xi)(x_1, \ldots , x_i) &~=~ f(\xi(x_1, \ldots , x_i)) \\
   (\pi f)(x_1, \ldots , x_i) &~=~ \pi(f(x_1), \ldots , f(x_i))
   \end{split}
   \]
for $x_1, \ldots , x_i \in R$.

\begin{lemma}
\label{lem:dd=0}
For $i = 1, 2$, we have $d^{i+1}_f d^i_f = 0$.
\end{lemma}

\begin{proof}
For $i = 1, 2$, we already know that $d^{i+1}d^i = 0$ in both $R$ and $S$.  Thus, to prove the Lemma, it remains to show that the third component of $d^{i+1}_f d^i_f(\xi; \pi; \varphi)$, namely, $f(d^i \xi) - (d^i\pi)f - d^i(f\xi - \pi f - d^{i-1}\varphi)$, is equal to $0$.  A straightforward inspection shows that $f(d^i\xi) = d^i(f\xi)$ and $(d^i\pi)f = d^i(\pi f)$.  This is enough, since $d^id^{i-1} = 0$ in $\Czinb^*(R,S)$ as well.
\end{proof}

In view of Lemma \ref{lem:dd=0}, we have a cochain complex $(\Czinb^*(f,f), d_f)$ $(1 \leq * \leq 4)$.  We call this the \emph{deformation complex of $f$}.  For $n = 2, 3$, define
   \[
   \Hzinb^n(f,f) := H^n(\Czinb^*(f,f), d_f).
   \]
The deformation complex and its cohomology groups will be used in the following sections to describe the behavior of $f$ with respect to algebraic deformations.

We will often use the shorter notation $\Czinb^*(f,f)$ to denote the deformation complex.  Since $\Czinb^0(R,S) = 0$ by definition, we will denote elements in $\Czinb^1(f,f)$ by pairs $(\xi,\pi) \in \Czinb^1(R,R) \times \Czinb^1(S,S)$.

%%=================================%%
%%                                 %%
%%  Deformation and infinitesimal  %%
%%                                 %%
%%=================================%%

\section{Deformation and infinitesimal}
\label{sec:def inf}

Fix a morphism $f \colon R \to S$ of dual Leibniz algebras.  Write $m_R$ (respectively, $m_S$) for the binary operation in $R$ (respectively, $S$).  In what follows, $t$ will denote an indeterminate.  As stated in the Introduction, we will use the results in \cite{bal3} concerning deformations of Zinbiel algebras freely.

%%%%%%%%%%%%%%%%%%%%%%%%
\subsection{Deformation}
\label{subsec:deformation}

A \emph{deformation of $f$} is a power series $\Theta_t = \sum_{i=0}^\infty \theta_i t^i$, in which $\theta_0 = (m_R; m_S; f)$ and each $\theta_i = (m_{R,i}; m_{S,i}; f_i) \in \Czinb^2(f,f)$, satisfying the following conditions:
   \begin{enumerate}
   \item $M_{R,t} = \sum_{i=0}^\infty m_{R,i}t^i$ is a deformation of $R$, i.e.
   \begin{equation}
   \label{eq:def}
   M_{R,t}(M_{R,t}(x,y),z) ~=~ M_{R,t}(x, M_{R,t}(y,z)) + M_{R,t}(x, M_{R,t}(z,y))
   \end{equation}
for all $x, y, z \in R$.
   \item $M_{S,t} = \sum_{i=0}^\infty m_{S,i}t^i$ is a deformation of $S$.
   \item $F_t = \sum_{i=0}^\infty f_it^i$ is a morphism of dual Leibniz algebras $(R \lbrack \lbrack t \rbrack \rbrack, M_{R,t}) \to (S \lbrack \lbrack t \rbrack \rbrack, M_{S,t})$.  Equivalently,
   \begin{equation}
   \label{eq:deformation}
   F_t(M_{R,t}(x,y)) ~=~ M_{S,t}(F_t(x), F_t(y))
   \end{equation}
   for all $x, y \in R$.
   \end{enumerate}
Instead of thinking of $\Theta_t$ as a power series in which each coefficient is a triple, we will sometimes treat it as a triple whose components are power series and write $\Theta_t = (M_{R,t}; M_{S,t}; F_t)$.

Note that, given $2$-cochains $m_{*,i} \in \Czinb^2(*,*)$ $(* = R,S)$ with $m_{*,0} = m_*$, $M_{*,t} = \sum_{i=0}^\infty m_{*,i}t^i$ is a deformation if and only if
   \begin{equation}
   \label{eq:def'}
   \sum_{l = 0}^n m_{*,l}(m_{*,n-l}(x, y), z) ~=~ \sum_{l=0}^n m_{*,l}(x, m_{*,n-l}(y, z) + m_{*,n-l}(z, y)) \tag{\ref{eq:def}$_n$}
   \end{equation}
for $x, y, z \in *$ and all $n \geq 0$.  Likewise, the condition \eqref{eq:deformation} is equivalent to
   \begin{equation}
   \label{eq:deformation'}
   \sum_{i=0}^n f_i m_{R,n-i}(x,y) ~=~ \sum_{i+j+k\,=\,n} m_{S,i}(f_j(x), f_k(y)) \tag{\ref{eq:deformation}$_n$}
   \end{equation}
for $x, y \in R$ and all $n \geq 0$.

A \emph{formal isomorphism of $f$} is a power series $\Phi_t = \sum_{i=0}^\infty (\phi_{R,i}; \phi_{S,i})t^i$, in which $(\phi_{R,0}; \phi_{S,0}) = (\Id_R; \Id_S)$ and each $(\phi_{R,i}; \phi_{S,i}) \in \Czinb^1(f,f) = \Czinb^1(R,R) \times \Czinb^1(S,S)$.  We will also denote $\Phi_t$ by $(\Phi_{R,t}; \Phi_{S,t})$.

Let $\Thetabar_t = \sum_{i=0}^\infty \thetabar_i t^i$, $\thetabar_i = (\mbar_{R,i}; \mbar_{S,i}; \fbar_i)$, be another deformation of $f$.  Then we say that the deformations $\Theta_t$ and $\Thetabar_t$ are \emph{equivalent} if and only if there exists a formal isomorphism $\Phi_t = (\Phi_{R,t}; \Phi_{S,t})$ of $f$ such that
   \begin{subequations}
   \label{eq:equivalence}
   \begin{align}
   \Mbar_{R,t} &:= \sum_{i=0}^\infty \mbar_{R,i} t^i = \Phi_{R,t} M_{R,t} \Phi_{R,t}^{-1}, \label{eq:equivalence R} \\
   \Mbar_{S,t} &:= \sum_{i=0}^\infty \mbar_{S,i} t^i = \Phi_{S,t} M_{S,t} \Phi_{S,t}^{-1}, \label{eq:equivalence S} \\
   \Fbar_t &:= \sum_{i=0}^\infty \fbar_it^i = \Phi_{S,t}F_t\Phi_{R,t}^{-1}.
   \end{align}
   \end{subequations}
In this case, we write $\Thetabar_t = \Phi_t \Theta_t \Phi_t^{-1}$.  Conversely, given only a deformation $\Theta_t$ and a formal isomorphism $\Phi_t$, one can define another deformation, $\Thetabar_t := \Phi_t \Theta_t \Phi_t^{-1}$, using \eqref{eq:equivalence}.  The resulting deformation is equivalent to $\Theta_t$ via $\Phi_t$.

The \emph{trivial deformation of $f$} is the deformation $\Theta_t = \theta_0 = (m_R; m_S; f)$.  The morphism $f$ is said to be \emph{rigid} if and only if every deformation of $f$ is equivalent to the trivial deformation.

\subsection{Infinitesimal}
\label{subsec:infinitesimal}

Let $\Theta_t = \sum_{i=0}^\infty \theta_i t^i$ be a deformation of $f$.  Its linear coefficient $\theta_1 = (m_{R,1}; m_{S,1}; f_1)$ is called the \emph{infinitesimal}.

\begin{thm}
\label{thm:inf}
The infinitesimal $\theta_1$ is a $2$-cocycle in $\Czinb^2(f,f)$.  More generally, if $\theta_i = 0$ for $i = 1, 2, \ldots, l$, then $\theta_{l+1}$ is a $2$-cocycle.  Moreover, the cohomology class of $\theta_1$ is well defined by the equivalence class of $\Theta_t$.
\end{thm}

\begin{proof}
Since $m_{R,1}$ is the infinitesimal of the deformation $M_{R,t}$ of $R$, it is a $2$-cocycle in $\Czinb^2(R,R)$.  The condition \eqref{eq:deformation}$_1$ is equivalent to $(fm_{R,1} - m_{S,1}f - d^1f_1) = 0$.  Therefore, $\theta_1$ is a $2$-cocycle.  The argument for the assertion concerning $\theta_{l+1}$ is similar, using \eqref{eq:deformation}$_{l+1}$ instead.

For the last assertion, let $\Phi_t$ be a formal isomorphism of $f$ and $\Thetabar_t = (\Mbar_{R,t}; \Mbar_{S,t}; \Fbar_t) = \Phi_t \Theta_t \Phi_t^{-1}$ be a deformation that is equivalent to $\Theta_t$.  Then $(m_{R,1} - \mbar_{R,1}) = d^1\phi_{R,1}$ and similar for $S$.  Computing modulo $t^2$, we have
   \[
   \begin{split}
   \Fbar_t(x)
   &\equiv (\Id_S + \phi_{S,1}t)(f + f_1 t)(x - \phi_{R,1}(x)t) \\
   &\equiv f(x) + (f_1(x) + \phi_{S,1}(f(x)) - f(\phi_{R,1}(x)))t.
   \end{split}
   \]
In particular, we have $f_1 - \fbar_1 = f\phi_{R,1} - \phi_{S,1}f$, which implies that $\theta_1 - \thetabar_1 = d^1_f(\phi_{R,1}; \phi_{S,1})$, a $2$-coboundary in $\Czinb^2(f,f)$.
\end{proof}

%%============%%
%%            %%
%%  Rigidity  %%
%%            %%
%%============%%

\section{Cohomological criterion for rigidity}
\label{sec:rigidity}

In this section, we obtain a simple cohomological criterion for the rigidity of the morphism $f$.  This will follow easily from the Theorem below.

\begin{thm}
\label{thm:rigidity}
Let $\Theta_t = \theta + \theta_l t^l + \theta_{l+1}t^{l+1} + \cdots $ be a deformation of $f$ for some $l \geq 1$ such that $\theta_l$ is a $2$-coboundary in $\Czinb^2(f,f)$.  Then there exists a formal isomorphism $\Phi_t$ of $f$ of the form $\Phi_t = (\Id_R + \phi_Rt^l; \Id_S + \phi_St^l)$
for some $\phi_R \in \Czinb^1(R,R)$ and $\phi_S \in \Czinb^1(S,S)$, such that the deformation $\Thetabar_t = \sum_{i=0}^\infty \thetabar_it^i := \Phi_t \Theta_t \Phi_t^{-1}$ satisfies $\thetabar_i = 0$ for $1 \leq i \leq l$.
\end{thm}

\begin{proof}
By hypothesis, there exists a $1$-cochain $(\phi_R; \phi_S) \in \Czinb^1(f,f)$ such that
   \[
   \begin{split}
   \theta_l
   &~=~ (m_{R,l}; m_{S,l}; f_l) \\
   &~=~ d^1(\phi_R; \phi_S) \\
   &~=~ (d^1 \phi_R; d^1 \phi_S; f\phi_R - \phi_Sf).
   \end{split}
   \]
Define the formal isomorphism $\Phi_t = (\Id_R + \phi_R t^l; \Id_S + \phi_S t^l)$ and the deformation $\Thetabar_t = (\Mbar_{R,t}; \Mbar_{S,t}; \Fbar_t) = \Phi_t \Theta_t \Phi_t^{-1}$.  Then$\Mbar_{R,t} = \sum_{i=0}^\infty \mbar_{R,i}t^i$ satisfies $\mbar_{R,i} = 0$ for $1 \leq i \leq l$.  The same holds with $S$ in place of $R$.  Computing modulo $t^{l+1}$, we have
   \[
   \begin{split}
   \Fbar_t(x)
   &~\equiv~ (\Id_S + \phi_S t^l)(f + f_l t^l)(x - \phi_R(x)t^l) \\
   &~\equiv~ f(x) + (f_l(x) + \phi_S(f(x)) - f(\phi_R(x)))t^l \\
   &~\equiv~ f(x).
   \end{split}
   \]
Therefore, $\fbar_i = 0$ for $1 \leq i \leq l$, from which the Theorem follows.
\end{proof}

The promised criterion for rigidity is now an immediate consequence of Theorem \ref{thm:inf} and Theorem \ref{thm:rigidity}.

\begin{cor}
\label{cor:rigidity}
If $\Hzinb^2(f,f) = 0$, then $f$ is rigid.
\end{cor}

%%==============%%
%%              %%
%%  Extensions  %%
%%              %%
%%==============%%

\section{Extending $2$-cocycles to deformations}
\label{sec:extending}

By Theorem \ref{thm:inf}, every infinitesimal is a $2$-cocycle in the deformation complex of $f: R \to S$.  It is, therefore, natural to ask if a given $2$-cocycle can be realized as the infinitesimal of some deformation.  In this section, we identify the obstructions for the existence of such a deformation.  We begin with some definitions.

\subsection{Deformations of finite order}
\label{subsec:finite}

Let $N$ denote a positive integer or $\infty$.  By a \emph{deformation of order $N$ of $f$}, we mean a power series $\Theta_t = \sum_{i=0}^N \theta_i t^i$, in which $\theta_0 = (m_R; m_S; f)$ and each $\theta_i = (m_{R,i}; m_{S,i}; f_i) \in \Czinb^2(f,f)$, satisfying (1) - (3) in Section \ref{subsec:deformation} modulo $t^{N+1}$ (where $t^{\infty} = 0$).  In other words, \eqref{eq:def'} ($* = R, S$) and \eqref{eq:deformation'} hold for $x, y, z \in R$ and $0 \leq n \leq N$.  Just as before, we will also write $\Theta_t = (M_{R,t}; M_{S,t}; F_t)$, where $F_t = \sum_{i=0}^N f_it^i$.

If $\Theta_t$ is a deformation of order $N < \infty$, then we say that it \emph{extends to order $N + 1$} if and only if there exists $\theta_{N+1} \in \Czinb^2(f,f)$ such that $\Thetatilde_t := \Theta_t + \theta_{N+1}t^{N+1}$ is a deformation of order $N + 1$.  In this case, $\Thetatilde_t$ is called an \emph{order $N + 1$ extension of $\Theta_t$}.

A deformation of order $\infty$ is just a deformation.  If $\theta_1 \in \Czinb^2(f,f)$ is a $2$-cochain, then $\Theta_t = \theta + \theta_1t$ is a deformation of order $1$ if and only if $\theta_1$ is a $2$-cocycle.  Thus, to answer the question discussed at the beginning of this section, it suffices to find the obstructions to extending a deformation of order $N$ to one of order $N + 1$ for $1 \leq N < \infty$.

\subsection{Obstructions}
\label{subsec:obstructions}

Now let $\Theta_t = \sum_{i=0}^N \theta_i t^i$ be a deformation of order $N < \infty$.  Consider the $3$-cochain
   \begin{equation}
   \label{eq:Ob def}
   \Ob_\Theta ~=~ (\Ob_R; \Ob_S; \Ob_f) \in \Czinb^3(f,f)
   \end{equation}
whose components are defined as follows.  For $* = R, S$ and $x, y, z \in R$, set
   \[
   \begin{split}
   \Ob_*(x,y,z)
   &~=~ \sum_{i=1}^N m_{*,i}(m_{*,N+1-i}(x,y),z) \\
   &\relphantom{} ~-~ \sum_{i=1}^N m_{*,i}(x, m_{*,N+1-i}(y,z) + m_{*,N+1-i}(z,y)) \\
   \Ob_f(x,y) &~=~ \sumprime m_{S,i}(f_j(x), f_k(y)) ~-~ \sum_{i=1}^N f_i m_{R,N+1-i}(x,y),
   \end{split}
   \]
where
   \[
   \sumprime ~=~
   \sum_{\substack{i+j\,=\,N+1 \\ k\,=\,0 \\ i,j\,>\,0}}
   ~+~ \sum_{\substack{i+k\,=\,N+1 \\ j\,=\,0 \\ i,k\,>\,0}}
   ~+~ \sum_{\substack{j+k\,=\,N+1 \\ i\,=\,0 \\ j,k\,>\,0}}
   ~+~ \sum_{\substack{i+j+k\,=\,N+1 \\ i,j,k\,>\,0}}.
   \]
The $3$-cochain $\Ob_\Theta$ is called the \emph{obstruction class of $\Theta_t$}.

\begin{lemma}
\label{lem:extending}
The obstruction class $\Ob_\Theta$ is a $3$-cocycle in $\Czinb^3(f,f)$.
\end{lemma}

The proof of this Lemma is a long but elementary computation, a sketch of which will be given in the next section.  The following Theorem does not depend on Lemma \ref{lem:extending}.

\begin{thm}
\label{thm:extending}
Let $\Theta_t$ be a deformation of order $N < \infty$ and $\theta_{N+1} = (m_{R,N+1}; m_{S,N+1}; f_{N+1})$ be a $2$-cochain in $\Czinb^2(f,f)$.  Then $\Thetatilde_t := \Theta_t + \theta_{N+1}t^{N+1}$ is an order $N + 1$ extension of $\Theta_t$ if and only if $\Ob_\Theta = d^2_f \theta_{N+1}$.
\end{thm}

\begin{proof}
Write $\Theta_t = (M_{R,t}; M_{S,t}; F_t)$ and  $\Thetatilde_t = (\Mtilde_{R,t}; \Mtilde_{S,t}; \Ftilde_t)$.  Then $\Thetatilde_t$ is a deformation of order $N + 1$ if and only if the following three conditions hold:
   \begin{enumerate}
   \item $\Mtilde_{R,t} = M_{R,t} + m_{R,N+1} t^{N+1}$ is an order $N + 1$ extension of $M_{R,t}$.
   \item $\Mtilde_{S,t} = M_{S,t} + m_{S,N+1} t^{N+1}$ is an order $N + 1$ extension of $M_{S,t}$.
   \item The condition \eqref{eq:deformation}$_{N+1}$ holds.
   \end{enumerate}
Condition (1) (respectively, (2)) is equivalent to $\Ob_R = d^2m_{R,N+1}$ (respectively, $\Ob_S = d^2m_{S,N+1}$).  Moreover, by a simple rearrangement of terms, it is easy to see that \eqref{eq:deformation}$_{N+1}$ is equivalent to $\Ob_f = (fm_{R,N+1} - m_{S,N+1}f - d^1f_{N+1})$.  This finishes the proof of the Theorem.
\end{proof}

In particular, assuming Lemma \ref{lem:extending} for the moment, $\Theta_t$ extends to order $N + 1$ if and only if the cohomology class of $\Ob_\Theta$ vanishes.  Starting with a $2$-cocycle and applying Lemma \ref{lem:extending} and Theorem \ref{thm:extending} repeatedly, we immediately obtain the following cohomological criterion for extending  $2$-cocycles to deformations.

\begin{cor}
\label{cor:extending}
If $\Hzinb^3(f,f)$ is trivial, then every $2$-cocycle in $\Czinb^2(f,f)$ is the infinitesimal of some deformation of $f$.
\end{cor}

We now turn to the proof of Lemma \ref{lem:extending}.

%%==============%%
%%              %%
%%  Main Lemma  %%
%%              %%
%%==============%%

\section{Sketch of the proof of Lemma \ref{lem:extending}}
\label{sec:Lemma}

In this final section, we give a sketch of the proof of Lemma \ref{lem:extending}.

We need to show that the obstruction class $\Ob_\Theta$ of $\Theta_t$ is a $3$-cocycle in the deformation complex of $f$.  We already know that $\Ob_R$ (respectively, $\Ob_S$) is a $3$-cocycle in $\Czinb^3(R,R)$ (respectively, $\Czinb^3(S,S)$).  Thus, it remains to show that
   \begin{equation}
   \label{eq:tobeshown}
   f\Ob_R \,-\, \Ob_S f ~=~ d^2\Ob_f
   \end{equation}
in $\Czinb^3(R,S)$.

To see this, let $x, y, z$ be elements in $R$.  Then the left-hand side of \eqref{eq:tobeshown} gives
   \begin{subequations}
   \allowdisplaybreaks
   \label{eq:LHS}
   \begin{align}
   (f \Ob_R & \,-\, \Ob_S f)(x,y,z) \notag \\
   &=~ \sum fm_{R,i}(m_{R,N+1-i}(x,y),z) \label{eq:A1} \\
   &\relphantom{} \relphantom{} -~ \sum fm_{R,i}(x, m_{R,N+1-i}(y,z) + m_{R,N+1-i}(z,y)) \label{eq:A2} \\
   &\relphantom{} \relphantom{} -~ \sum m_{S,i}(m_{S,N+1-i}(f(x),f(y)), f(z)) \label{eq:B1} \\
   &\relphantom{} \relphantom{} +~ \sum m_{S,i}(f(x), m_{S,N+1-i}(f(y),f(z))) \label{eq:B21} \\
   &\relphantom{} \relphantom{} +~ \sum m_{S,i}(f(x), m_{S,N+1-i}(f(z),f(y))) \label{eq:B22},
   \end{align}
   \end{subequations}
where each sum is $\sum_{i=1}^N$.  On the other hand, the right-hand side of \eqref{eq:tobeshown} gives
   \begin{subequations}
   \allowdisplaybreaks
   \label{eq:RHS}
   \begin{align}
   (d^2 \Ob_f)(x,y,&z) \notag \\
   &=~ f(x) \cdot \Ob_f(y,z) \,+\, f(x) \cdot \Ob_f(z,y) \label{eq:1} \\
   &\relphantom{} \relphantom{} -~ \Ob_f(x \cdot y, z) \label{eq:2} \\
   &\relphantom{} \relphantom{} +~ \Ob_f(x, y\cdot z \,+\, z \cdot y) \label{eq:3} \\
   &\relphantom{} \relphantom{} -~ \Ob_f(x,y) \cdot f(z). \label{eq:4}
   \end{align}
   \end{subequations}
The rest of the proof is an elementary computation that involves expanding and rewriting the terms \eqref{eq:1}, \ldots , \eqref{eq:4}, using \eqref{eq:def'} and \eqref{eq:deformation'} for various vales of $n$ repeatedly, and observing that their sum is the sum of \eqref{eq:A1}, \ldots , \eqref{eq:B22}.

For example, consider \eqref{eq:3}.  By the definition of $\Ob_f$ and the linearity of $f_k$ ($k \geq 0$), we have
   \begin{subequations}
   \label{eq':3}
   \allowdisplaybreaks
   \begin{align}
   \Ob_f(x,\, & y\cdot z + z \cdot y) \notag \\
   &=\, \sum_{\substack{i+j\,=\,N+1 \\ i,j\,>\,0}} m_{S,i}(f_j(x), f(y \cdot z) \,+\, f(z \cdot y)) \label{eq:3a} \\
   &\relphantom{} +\, \sum_{\substack{i+k\,=\,N+1 \\ i,k\,>\,0}} m_{S,i}(f(x), f_k(y \cdot z) \,+\, f_k(z \cdot y)) \label{eq:3b} \\
   &\relphantom{} +\, \sum_{\substack{j+k\,=\,N+1 \\ j,k\,>\,0}} f_j(x) \cdot (f_k(y \cdot z) \,+\, f_k(z \cdot y)) \label{eq:3c} \\
   &\relphantom{} +\, \sum_{\substack{i+j+k\,=\,N+1 \\ i,j,k \,>\,0}} m_{S,i}(f_j(x), f_k(y \cdot z) \,+\, f_k(z \cdot y)) \label{eq:3d} \\
   &\relphantom{} -\, \sum_{i=1}^N f_i m_{R,N+1-i}(x, y \cdot z \,+\, z \cdot y). \label{eq:3e}
   \end{align}
   \end{subequations}
The term \eqref{eq:3d} needs to be expanded.  Using \eqref{eq:deformation}$_k$, we have that
   \[
   \begin{split}
   f_k(y \cdot z) & \,+\, f_k(z \cdot y) \\
   &=\, \sum_{a+b+c\,=\,k} \biglpren m_{S,a}(f_b(y),f_c(z)) \,+\, m_{S,a}(f_b(z),f_c(y)) \bigrpren \\
   &\relphantom{} -\, \sum_{l=0}^{k-1} \left( f_l m_{R,k-l}(y,z) \,+\, f_l m_{R,k-l}(z,y) \right).
   \end{split}
   \]
Putting this back into \eqref{eq:3d}, we have
   \begin{subequations}
   \label{eq':3d}
   \allowdisplaybreaks
   \begin{align}
   &\sum_{\substack{i+j+k\,=\,N+1 \\ i,j,k \,>\,0}} m_{S,i}(f_j(x), f_k(y \cdot z) \,+\, f_k(z \cdot y)) \notag \\
   &=\, \sum_{\substack{i+j+k\,=\,N+1 \\ i,j,k \,>\,0}} \sum_{a+b+c\,=\,k} m_{S,i}(f_j(x), m_{S,a}(f_b(y),f_c(z)) \,+\, m_{S,a}(f_b(z),f_c(y))) \label{eq:3d1} \\
   &\relphantom{} -\, \sum_{\substack{i+j+k\,=\,N+1 \\ i,j,k \,>\,0}} \sum_{l=0}^{k-1} m_{S,i}(f_j(x), f_l m_{R,k-l}(y,z) \,+\, f_l m_{R,k-l}(z,y)). \label{eq:3d2}
   \end{align}
   \end{subequations}
Therefore, \eqref{eq:3} can be rewritten as the sum of \eqref{eq:3a}, \eqref{eq:3b}, \eqref{eq:3c}, \eqref{eq:3e}, \eqref{eq:3d1}, and \eqref{eq:3d2}.

The same type of argument can be applied to \eqref{eq:2}.  Once this is done, one writes down the terms \eqref{eq:1} and \eqref{eq:4} using the definition of $\Ob_f$.  The required condition \eqref{eq:tobeshown} will then follow immediately.

%%==============%%
%%              %%
%%  References  %%
%%              %%
%%==============%%

\end{document}